\documentclass[10pt,reqno]{amsart}

\usepackage[a4paper,asymmetric]{geometry}
\usepackage{latexsym, amsmath, amssymb, esint}
\usepackage{bm}
\usepackage{graphics,epsfig,graphicx,graphics}
\usepackage{color, xcolor}
\usepackage{tikz}
\usepackage[hyperpageref]{backref}
\usetikzlibrary{patterns}
\usepackage[english]{babel}
\usepackage[colorlinks=true,urlcolor=blue, citecolor=red,linkcolor=blue,
linktocpage,pdfpagelabels, bookmarksnumbered,bookmarksopen]{hyperref}
\newtheorem{theorem}{Theorem}[section]

\newtheorem{prop}[theorem]{Proposition}

\theoremstyle{definition}

\newcommand{\ep}{\varepsilon}
\newcommand{\R}{\mathbb{R}}
\newcommand{\rn}{\R^N}
\newcommand{\N}{\mathbb{N}}

\newcommand{\bd}{\partial}

\newcommand{\bho}{B_{H_0}}
\newcommand{\bhor}[1]{B_{H_0}({#1})}
\newcommand{\bhorx}[2]{B_{H_0}({#2},{#1})}

\newcommand{\CC}{\textsf{C }}

\newcommand{\Hcap}{{\rm Cap_{H,p}}}
\newcommand{\hnorm}{{\mathcal{J}_p}}

\newcommand{\la}{\lambda}
\newcommand{\uom}{u_\Omega}

\newcommand{\qq}{{\textbf{\textsf q}}}%{q_c}%{\overline{q}}
\newcommand{\ts}{\textsf{t}} %{\overline{t}}
\newcommand{\ls}{\textsf{s}} %{\overline{s}}

\numberwithin{equation}{section}
\begin{document}
	\title[Overdetermined problems for the anisotropic capacity]{Some overdetermined problems related to the anisotropic capacity}
	
	%\author[C. Bianchini]{Chiara Bianchini}
	%\author[G. Ciraolo]{Giulio Ciraolo}
	
	\author[C. Bianchini]{Chiara Bianchini}
	\author[G. Ciraolo]{Giulio Ciraolo}
	\author[P. Salani]{Paolo Salani}

	\address{C. Bianchini, Dipartimento di Matematica e Informatica ``U. Dini'', Universit\`a degli Studi di Firenze, Viale Morgagni 67/A, 50134 Firenze - Italy}
	\email{chiara.bianchini@unifi.it}

	\address{G. Ciraolo,  Dipartimento di Matematica e Informatica, Universit\`a degli Studi di Palermo, Via Archirafi 34, 90123, Palermo - Italy}
	\email{giulio.ciraolo@unipa.it}

	\address{P. Salani,  Dipartimento di Matematica e Informatica ``U. Dini'', Universit\`a degli Studi di Firenze, Viale Morgagni 67/A, 50134 Firenze - Italy}
	\email{paolo.salani@unifi.it}

	\begin{abstract}
	We characterize the Wulff shape of an anisotropic norm in terms of solutions to overdetermined problems for the Finsler $p$-capacity of a convex set $\Omega \subset \rn$, with $1<p<N$. In particular we show that if the Finsler $p$-capacitary potential $u$ associated to $\Omega$ has two homothetic level sets then $\Omega$ is Wulff shape. Moreover, we show that the concavity exponent of $u$ is $\qq=-(p-1)/(N-p)$ if and only if $\Omega$ is Wulff shape.
	\end{abstract}
	
	\maketitle

	

	%%%%%%%%%%%%%%%%%%%%%%%%%%%%%%%%%%%%%%%%%%%%%%%%%
	%%%%%%%%%%%%%%%%%%%%%%%%%%%%%%%%%%%%%%%%%%%%%%%%%

%%%%%%%%%%%%%%%%%%%%%%%%%%%%%%%%%%%%%%%%%%%%%%%%%%%%%%%%%%%%%%%%%%%%

%%%%%%%%%%%%%%%%%%%%%%%%%%%%%%%%%%%%%%%%%%%%%%%%%
% % % % % % % % % % % % % % % % % % % % % % % % % % % % % % % % % % % % % % % % % % % % % % % % % % % % % % %
\section{Introduction}
% % % % % % % % % % % % % % % % % % % % % % % % % % % % % % % % % % % % % % % % % % % % % % % % % % % % % % % %

%
%
The aim of this paper is to study some unconventional overdetermined problems for the Finsler $p$-capacity of a bounded convex set $\Omega$ associated to a norm $H$ of $\rn$, $N\geq 3$. 

Given a bounded convex domain $\Omega \subset \rn$, the $p$-capacity of $\Omega$ is defined by
\begin{equation*}
{\rm Cap_p}(\Omega)=\inf \left\{  \frac 1p \int_{\rn} |D\varphi|^p \, dx,\ \varphi\in C^{\infty}_0(\rn), \varphi(x)\ge 1 \text{ for } x\in\Omega \right\} 
\end{equation*}
with $1<p<N$. When the Euclidean norm $|\cdot|$ is replaced by a more general norm $H(\cdot)$, one can consider the so called \emph{Finsler} $p$-capacity $\Hcap{(\Omega)}$, which is defined by
\begin{equation}\label{pcap}
\Hcap(\Omega)=\inf \left\{  \frac 1p\int_{\rn} H^p(D\varphi)\, dx,\ \varphi\in C^{\infty}_0(\rn), \varphi(x)\ge 1 \text{ for } x\in\Omega \right\} \,,
\end{equation}
for $1<p<N$.
Under suitable assumptions on the norm $H$ and on the set $\Omega$, the above infimum is attained and
$$
\Hcap(\Omega)=\frac 1p\int_{\rn} H^p(D\uom)\, dx\,,
$$
where $\uom$ is the solution of the Finsler $p$-capacity problem
\begin{equation}\label{ext-pbu}
\begin{cases}
\Delta^H_p u =0\qquad&\text{in }\rn\setminus\overline{\Omega},\\
u=1\qquad&\text{on }\bd\Omega \,,\\
u\to 0\qquad&\text{as } H(x)\to +\infty \,.
\end{cases}
\end{equation}
Here $\Delta^H_p$ denotes the Finsler $p$-Laplace operator, i.e. $\Delta^H_p u=\text{div}(H^{p-1}(Du)\nabla H(Du))$.
The function $u_{\Omega}$ is named  \emph{(Finsler) $p$-capacitary potential} of $\Omega$.

When $\Omega$ is Wulff shape, i.e. it is a sublevel set of the dual norm $H_0$
$$
\Omega=\bhor{r} = \{x \in \rn:\ H_0(x - \bar x)<r \}
$$  
(see Section \ref{secNotations} for definitions), the solution to \eqref{ext-pbu} can be explicitly computed and it is given by 
\begin{equation}\label{potenziale-palla}
v_r(x)=\Big(\frac{H_0(x-\bar{x})}r\Big)^{\frac1 \qq}\,,
\end{equation}
with
\begin{equation}\label{q}
\qq = -\frac{p-1}{N-p} \,.
\end{equation}
It is straightforward to verify that the potential $v_r$ in \eqref{potenziale-palla} enjoys the following properties:
\begin{itemize}
\item[(i)] the function $v_r^{\qq}$ is convex, i.e. $v_r$ is {\em $\qq$-concave};
\item[(ii)] the superlevel sets of $v_r$ are homothetic sets and they are Wulff shapes;
\item[(iii)] $H(D v_r)$ is constant on the level sets of $v_r$.  
\end{itemize}

\noindent The aim of this paper is to show that each of the properties (i)-(iii) characterizes the Wulff shape under some regularity assumptions on the norm $H$ and on $\Omega$. In particular, we assume that $H \in \hnorm$ where 
\begin{equation}\label{ip}
\hnorm=\{ H\in C^{2}_+ (\rn\setminus\{0\}) \,, H^p\in C^{2,1}(\rn\setminus\{0\}) \ \} \,.
\end{equation}

Our first main result is related to property (i) and it is about concavity properties of the solution to \eqref{ext-pbu}. We recall that a nonnegative function $v$ with convex support is \emph{$\alpha$-concave}, for some $\alpha \in [-\infty, + \infty]$, if  
\begin{itemize}
\item $v$ is a positive constant in its support set, in case $\alpha=+\infty$;
\item $v^\alpha$ is concave, in case $\alpha >0$;
\item $\log v$ is concave, in case $\alpha =0$ (and $v$ is called \emph{log-concave});
\item $v^\alpha$ is convex, in case $\alpha <0$;
\item all its super level sets $\{v>t\}$ are convex, in case $\alpha = - \infty$ (and $v$ is called \emph{quasi-concave}).
\end{itemize}
Notice that if $v$ is $\alpha$-concave for some $\alpha>-\infty$, then it is $\beta$-concave for every $\beta\in[-\infty,\alpha]$. Then quasi-concavity is the weakest among concavity properties.

Concavity properties of solutions to elliptic and parabolic equations are a popular field of investigation. Classical results in this framework are for instance the log-concavity of the first Dirichlet eigenfunction  of the Laplacian (see \cite{BL}), the preservation of concavity by the heat flow (see again \cite{BL}), the $\frac 12$-concavity (i.e. the concavity of the square root) of the torsion function (see \cite{ML, Bo, CF, Ke}), the quasi-concavity of the Newton potential and of the $p$-capacitary potential (see \cite{G, L}). 
The latter results are especially related to the situation we consider in this paper.
Indeed, it is proved in \cite{AGHLV} (and it can be also obtained with the methods of  \cite{BLS}) that when $\Omega$ is a convex domain, its $p$-capacitary potential $u_{\Omega}$ is a \emph{quasi-concave function}, i.e. all its superlevel sets are convex. 
Moreover, as we have seen, quasi-concavity is the weakest property in this context and one may expect and ask more than this. 
%Indeed if a function is $\alpha$-concave, then it is $\beta$-concave too, for every $\beta<\alpha$.
Then, following \cite{Ke, Sa, HNST}, it is natural to define the \emph{concavity exponent} associated to the solution to \eqref{ext-pbu} as
\begin{equation} \label{alpha}
\alpha (\Omega,p) = \sup \{\beta \leq  1 \, : \ u_\Omega \, \text{ is $\beta$-concave } \}\,.
\end{equation}
In the Euclidean case, it was proved in \cite{Sa} that the concavity exponent attains its maximum when $\Omega$ is a ball (and only in this case). %and an analogous question for the torsion function is treated in \cite{HNST}. 
In the following theorem, we characterize the Wulff shape in terms of property (i) above. More precisely, we generalize the results of \cite{Sa} to the anisotropic setting and we prove that the exponent $\qq$ characterizes the Wulff shape.

\begin{theorem} \label{main1}
	Let $H$ be a norm of $\rn$ in the class (\ref{ip}) and let $\Omega$ be a bounded convex domain of $\rn$ of class $C^2$. Then
	$$
	\alpha(\Omega,p) \le \qq, 
	$$
	with $\qq$ given by \eqref{q}, and equality holds if and only if $\Omega$ is Wulff shape.
\end{theorem}

The proof of Theorem \ref{main1} is based on the Brunn-Minkowski inequality for Finsler $p$-capacity, recently proved in \cite{AGHLV} (and here recalled in Proposition \ref{BMprop}),  and upon the fact $u$ can have a level set homothetic to $\Omega$ if and only if $\Omega$ is a ball. Clearly this property is related to property (ii) above. And indeed the characterization of the Wulff shape is achieved whenever \emph{just two} superlevel sets of $u_\Omega$ are homothetic, as expressed in the following theorem.

\begin{theorem}\label{thm-livelloomot}
Let $H$ be a norm of $\rn$ in the class (\ref{ip}). Let $\Omega \subset \mathbb{R}^N$, $N\geq 3$, be a bounded convex domain with boundary of class $C^2_+$. If there exists a solution to \eqref{ext-pbu} having two homothetic superlevel sets, then $\Omega$ is Wulff shape.	
\end{theorem}

The Euclidean counterpart of Theorem \ref{thm-livelloomot} was proved in \cite{Sa}. The proof of Theorem \ref{thm-livelloomot} in the anisotropic setting passes through the following theorem, which is related to property (iii).

\begin{theorem} \label{thm_overdet}
	Let $\Omega \subset \rn$ be a convex domain containing the origin and with boundary of class $C^{2,\alpha}$. Let $H \in \hnorm$ and let $R>0$ be such that $\overline{\Omega} \subset \bhor{R}$.   There exists a solution to   
	\begin{equation}\label{pbB_overdet}
	\begin{cases}
	\Delta_p^H u =0& \text{ in } \bhor{R} \setminus\overline{\Omega}\\
	u=1& \text{ on }\bd \Omega\\
	u = 0& \text{ on }\bd \bhor{R}  \\
	H(Du) = C & \text{ on }\bd \bhor{R} \,,
	\end{cases}
	\end{equation}
	for some constant $C>0$ if and only if $\Omega=\bhor{r}$, with 
	\begin{equation} \label{r_e_C}
	r\;C=\frac{N-p}{p-1} \,.
	\end{equation} 
\end{theorem}

Since two boundary conditions (Dirichlet and Neumann) are imposed on a prescribed part of the boundary, Theorem \ref{thm_overdet} clearly falls in the realm of {\em overdetermined problems}: since the domain $\Omega$ is not prescribed,  the unknown of the problem is in fact the couple $(\Omega, u)$, and by imposing that $u$ has some peculiar property (which is not commonly shared by all the solution of the involved PDE), one ask whether this is sufficient to uniquely determine the domain $\Omega$. In this sense, also Theorem \ref{main1} and Theorem \ref{thm-livelloomot} can be considered as overdetermined problems, since we ask for a solution of a Dirichlet problem satisfying some extra special condition ($\qq$-concavity or homothety of level sets, respectively).

\medskip

The paper is organized as follows. In Section \ref{secNotations} we introduce some notation and basic properties of Finsler norms; then we recall some known fact about the Finsler capacity $\Hcap$ of a convex set and, in particular, the Brunn-Minkowski inequality from \cite{AGHLV}.  Theorems \ref{thm_overdet}, \ref{thm-livelloomot} and \ref{main1} are proved in Sections \ref{proof3}, \ref{proof2} and \ref{proof1}, respectively.

\medskip

{\bf{Acknowledgements.}} The authors have been partially supported by GNAMPA of INdAM and by a PRIN Project of Italian MIUR.

% % % % % % % % % % % % % % % % % % % % % % % % % % % % % % % % % % % % % % % % % %
% % % % % % % % % % % % % % % % % % % % % % % % % % % % % % % % % % % % % % % % % %
%%%%%%%%%%%%%%%%%%%%%%%%%%%%%%%%%%%%%%%%%%%%%%%%%
\section{Notations}\label{secNotations}

\subsection{Norms of $\rn$}

We consider the space $\rn$ endowed with a generic norm $H: \rn \to \R$ such that:
\begin{itemize}
	\item[(i)] $H$ is convex;
	\item[(ii)] $H(\xi) \geq 0$ for $\xi \in \rn$ and $H(\xi)=0$ if and only if $\xi=0$;
	\item[(iii)] $H(t\xi) = |t| H(\xi)$ for $\xi \in \rn$ and $t\in \R$.
\end{itemize}

Then we identify the dual space of $\rn$ with $\rn$ itself via the scalar product $\langle \cdot;\cdot\rangle$.
Accordingly the space $\rn$ turns out to be endowed also with the dual norm $H_0$ given by
\begin{equation}\label{defH0}
H_0(x)=\sup_{\xi\neq 0}\frac{\langle {x};{\xi}\rangle}{H(\xi)}\quad\text{ for
}x\in\rn\,.
\end{equation}

We denote by $\bhor{r}$ the anisotropic ball centered at $O$ with radius $r$ in the norm $H_0$, i.e.
\begin{equation*}
\bhor{r} = \{x \in \rn:\ H_0(x) < r\}.
\end{equation*}
Analogously, we define 
\begin{equation*}
B_H(r) = \{\xi \in \rn:\ H(\xi) < r\}.
\end{equation*}
The sets $\bhor{r}$ and $B_H(r)$ are called \emph{Wulff shape} of $H$ and $H_0$, respectively; in the special case $r=1$ they are indicated by $\bho, B_H$, respectively. Notice that, in the language of the theory of convex bodies, $H$ is {\em the support function} of $\bho$ and $H_0$ is in turn the support function of $B_H$.

For a regular convex domain $\Omega$ the Finsler perimeter is defined by
$$
P_H(\bd\Omega)=\int_{\bd\Omega} H(\nu)\; d\sigma,
$$
where $\nu$ is the outer unit normal to $\bd\Omega$.

%For a domain $\Omega$ we indicate by $\conv(\Omega)$ its convex hull. 
\subsection{Finsler capacity}
For a bounded convex domain $\Omega$ in $\rn$ its \emph{Finsler p-capacity}, denoted by $\Hcap{(\Omega)}$, is defined as follows:
\begin{equation*} 
\Hcap(\Omega)=\inf \left\{  \frac 1p\int_{\rn} H^p(D\varphi)\, dx,\ \varphi\in C^{\infty}_0(\rn), \varphi(x)\ge 1 \text{ for } x\in\Omega \right\} \,,
\end{equation*}
for $N\geq 3$ and $1<p<N$.
If $H$ is a norm in the class (\ref{ip}), the integral operator is strictly convex and hence \eqref{pcap} admits a unique solution $u_{\Omega}$, which satisfies
\begin{equation*} 
\begin{cases}
\Delta^H_p u =0\qquad&\text{in }\rn\setminus\overline{\Omega},\\
u=1\qquad&\text{on }\bd\Omega \,,\\
u\to 0\qquad&\text{as } H(x)\to +\infty \,.
\end{cases}
\end{equation*}
The function $u_{\Omega}$ is called the \emph{Finsler p-capacitary potential} of $\Omega$.
As already noticed when $\Omega$ is a convex set the potential $u_{\Omega}$ is at least quasi-concave, that is its superlevel sets are convex sets (see Lemma 4.4 \cite{AGHLV}). 

In the special case $\Omega=\bhor{r}$ the capacitary potential is easily computed and is given by \eqref{potenziale-palla},
but this is not possible for general convex domain. However, when $\Omega$ is a convex set, asymptotic estimates for $u_{\Omega}$ are known. In particular, it has  recently been proved in \cite{AGHLV} the following:
\begin{equation}\label{limite}
\lim_{|x|\to \infty} u_{\Omega}(x)H_0(x)^\frac{N-p}{p-1} = \CC \Hcap^{\frac 1{p-1}}(\Omega),
\end{equation}
where $\CC=(N-2) P_H^{\frac 1{p-1}}(\bd\bho)$. Moreover, one can prove that there exists a positive constant $\gamma$ such that 
\begin{equation} \label{asy}
\gamma^{-1} H(x)^{\frac 1\qq -1} \leq H(Du(x)) \leq \gamma H(x)^{\frac 1\qq -1} \,,
\end{equation}
(see \cite{BCS}, \cite{BiCi}).

The $p$-capacity operator satisfies a Brunn-Minkowski inequality. In the Euclidean setting, this was proved in \cite{CS}, such result has been recently extended to quite general operators in divergence form in \cite{AGHLV}. Here, we recall the following from \cite{AGHLV}.
\begin{prop}[\cite{AGHLV}]\label{BMprop}
	Let $K,D$ be compact convex sets in $\rn$ satisfying 
	$$
	\Hcap(K),\Hcap(D)>0 \,.
	$$ 
	For $1<p<N$ and $\la\in[0,1]$ it holds
	\begin{equation}\label{BM}
	\Hcap^{\frac 1{N-p}}((1-\la)K+\la D) \ge (1-\la)\Hcap^{\frac 1{N-p}}(K)+\la\Hcap^{\frac 1{N-p}}(D),
	\end{equation}
	and equality holds if and only if $K$ and $D$ are homothetic sets.
\end{prop}

% % % % % % % % % % % % % % % % % % % % % % % % % % % % % % % % % % % % % % % % % % % % % % % % % % % % % % %

\section{Proof of Theorem \ref{thm_overdet}} \label{proof3}

	Let 
	$$
	v(x) = \frac{H_0(x)^{\frac 1 \qq} - R^{\frac 1 \qq}}{r^{\frac 1 \qq} - R^{\frac 1 \qq}} \,, \quad x \in \mathbb{R}^n\setminus \{O\} \,,
	$$
	with $r$ and $\qq$ given by \eqref{r_e_C} and \eqref{q}, respectively.
	When $\Omega=\bhor{r}$ is Wulff shape of radius $r$, a direct check shows that $v$ is the solution to \eqref{pbB_overdet}. 
	
	Now we prove the reverse assertion. Let 
	$$
	\tilde u (x) = 
	\begin{cases}
	u(x) & \text{if } x \in \overline{\bhor{R}} \setminus \Omega \,, \\
	v(x) & \text{if } x \in \mathbb{R}^n \setminus \overline{\bhor{R}} \,.
	\end{cases}
	$$
	We notice that $\tilde u \in C^1( \mathbb{R}^n \setminus \Omega) $ and it satisfies $\Delta_p^H \tilde u = 0$ in  $ \mathbb{R}^n \setminus \overline \Omega$ (which follows from the weak formulation of the equation).
	
	Fix any $t>1$ and set
	$E =  \bhor{tR} \setminus \overline{\Omega}$.
	For $\tau \in [0,1]$, we define 
	$$
	u_\tau = \tau \tilde u + (1-\tau) v 
	$$ 
	in $\overline E$; notice that $u_1 = \tilde u$ and $u_0 = v$.  
	
	The function $\tilde u - v$ satisfies an elliptic equation. Indeed, for any $\phi \in C_0^1 (E)$ we have
	\begin{equation*}
	\begin{split}
	0 & = \int_E H(D u_1 )^{p-1} \langle \nabla H(D u_1) ; D \phi\rangle \, dx - \int_E H(D u_0 )^{p-1} \langle \nabla H(D u_0) ; D \phi \rangle\, dx \\
	& = \int_E \langle\left( \int_0^1 \frac{d}{d\tau}  \left( H(D u_\tau )^{p-1} \nabla H(D u_\tau) \right)   \, d\tau \right) ; D \phi\rangle \, dx\\
	& = \int_E {\bf A}(x) \langle D (\tilde u - v) ; D \phi \rangle \, dx \,,
	\end{split}
	\end{equation*}
	where
	${\bf A}(x)=a_{ij}(x)$ is given by 
	\begin{equation*}
	\begin{split}
	a_{ij} (x) & = \int_0^1 \left( (p-1) H(Du_\tau)^{p-2} H_{\xi_j}(Du_\tau) H_{\xi_i}(Du_\tau) + H(Du_\tau)^{p-1} H_{\xi_i \xi_j}(Du_\tau) \right) d\tau \\
	& = \frac{1}{p}  \int_0^1 (\nabla^2 H^p (Du_\tau))_{ij} d\tau \,.
	\end{split}
	\end{equation*}
	From \eqref{asy} we have that
	\begin{equation} \label{altro}
	A \leq H(Du) \leq B
	\end{equation}
	in $E$ for some constant $A,B>0$.
	Notice that, since the super level sets of $u$ are convex sets (see \cite{AGHLV}) and  those of $v$ are Wulff shapes centered at the origin,
	we have
	\begin{equation}\label{cacca}
	\langle Du; \frac{x}{|x|} \rangle \ge \ep>0, \qquad \ep\leq\langle Dv; \frac{x}{|x|} \rangle \leq M,
	\end{equation}
	in $E$ for some positive constants $\ep$ and $M$. Hence, we can find $\tau_0 \in (0,1)$ such that
	$$(1-\tau_0)|Dv| \leq \tau_0 |Du|/2\,,%\,\,\,\text{ for any }\tau \in [\tau_0,1]\,,
	$$ 
	which implies that $|Du_\tau| \geq \tau_0 |Du|/2$  for every $\tau \in [\tau_0,1]$. From \eqref{cacca} and \eqref{altro} we obtain
	$$
	|Du_\tau| \geq \min\left((1-\tau_0)\ep,\, \frac{\tau_0}{2} A \right) > 0 \,\,\,\text{ in }E\,,
	$$
	which finally gives that  $a \le|Du_\tau|\le b$ in $E$
	for some constants $a,b>0$ . Such estimates imply that the operator
	$$
	Lw={\rm div} \left( {\bf A(x)} Dw \right)
	$$
	is uniformly elliptic. Moreover, since $H^p \in C^{2,1}(\mathbb{R}^n \setminus \{O\})$, we also have that $a_{ij}$ are locally Lipschitz. Hence, $L$ satisfies the assumptions of Theorem 1.1 in \cite{GL2} (see also \cite{GL1}) and we have the analytic continuation for $\tilde u - v$ in $E$, whence $\tilde u - v \equiv 0$ in $E$ which implies that $u \equiv v$ in $\overline{\bhor{R}}\setminus\Omega$ and we conclude.

%%%%%%%%%%%%%%%%%%%%%%%%%%%%%%%%%%%%%%%%%%%%%%%%%%%%%%%%%%%%%%%%%%

\section{Proof of Theorem \ref{thm-livelloomot}} \label{proof2}
% % % % % % % % % % % % % % % % % % % % % % % % % % % % % % % % % % % % % % % % % % % % % % % % % % % % % % % %

% % % % % % % % % %

	For any $t\in(0,1)$ we set $U(t)=\{ x\in\rn\ :\ u(x)\ge t \}$ and let $u_{U(t)}$ be the $p$-capacitary potential of $U(t)$. Hence 	\begin{equation}\label{ut}
	u_{U(t)}(x)= \frac 1t u(x),
	\end{equation}
	for every $x\in \rn\setminus  U(t)$, where $u$ is the Finsler $p$-capacitary potential of $\Omega$.
	
	Let $\ts, \ls \in (0,1)$, with $\ts < \ls$, be the levels of $u$ such that $U(\ts), U(\ls)$ are homothetic, that is: there exist $\xi\in\rn$ and $\rho>1$ such that $U(\ts)=\rho U(\ls) +\xi$. Up to a translation we can assume $\xi=0$ and  hence
	\begin{equation}\label{uUt}
	u_{U(\ts)}(x)=u_{U(\ls)}\Big( \frac{x}{\rho} \Big).
	\end{equation}
	
	\noindent \emph{Step 1:} $\rho^{-\frac{N-p}{p-1}}=\ts / \ls$.

	From (\ref{limite}), (\ref{ut}) and \eqref{uUt}, we have
	$$
	\CC \Hcap(\Omega)^{\frac{1}{p-1}}=\lim_{|x|\to\infty} \ts\ u_{U(\ls)}\Big(\frac{x}{\rho}\Big)H_0^{\frac{N-p}{p-1}}(x) \,.
	$$
	By using again (\ref{limite}) and (\ref{ut}), and from the homogeneity of $H_0$, we find
	$$
	\CC \Hcap(\Omega)^{\frac{1}{p-1}} = \frac{\ts}{\ls} \CC\Hcap(\Omega)^{\frac 1{p-1}} \rho^{\frac{N-p}{p-1}} % \lim_{|x|\to\infty} \frac{H_0^{\frac{N-p}{p-1}}(x)}{H_0^{\frac{N-p}{p-1}}(x-\xi)},
	$$
	which implies that $\rho^{-\frac{N-p}{p-1}}=\ts / \ls$.
	
	\medskip

	\noindent\emph{Step 2:} \emph{Let $r_k=\ts^k \;\ls^{1-k} $ for $k\ge 0$. Then $U(r_0)=U(\ls)$ and
		$U(r_k)=\rho^k U(\ls)$ for $k\in\N \,.$
	}
	
	Indeed notice that for every $z<\ts$ the set $U(z)$ is homothetic to $U(z \frac \ls{\ts})$ since by (\ref{ut}), (\ref{uUt}) we have
	$$
	U(z)=\{u(x)\ge z\} = \{ u_{U(\ls)}(\frac{x}{\rho}) \ge \frac{z}{\ts} \}= \{ u(\frac{x}{\rho}) \ge z \frac{\ls}{\ts} \},
	$$
	that is $U(z)=\rho U(z\frac{\ls}{\ts})$. 
	Hence, recalling that $U(\ls)=U(r_0)$ and that $r_k=\frac{\ts}{\ls} r_{k-1}$, we obtain  
	$$
	U(r_k)= \rho^k U(\ls) %+ \xi \frac{1-\rho^k}{1-\rho} \,,
	$$
	 for every $k\ge 0$.
	
	\medskip

	\noindent \emph{Step 3:} \emph{$U(\ls)$ is Wulff shape}.
	
	Let $x,y\in\bd U(\ls)$ and define 
	\begin{eqnarray*}
		x_k= \rho^k x\,,\\ %+ \xi \frac{1-\rho^k}{1-\rho},\\
		y_k=\rho^k y \,.%+ \xi \frac{1-\rho^k}{1-\rho}.
	\end{eqnarray*}
	Notice that 
	\begin{equation} \label{A}
	\lim_{k\to\infty}|x_k| = \lim_{k\to\infty}|y_k| = +\infty \,.
	\end{equation}
	From \emph{Step 2} the points $x_k,y_k$ belong to $\bd U(r_k)$, so that $u(x_k)=u(y_k)=r_k$. From \eqref{A} and  (\ref{limite}) we obtain
	$$
	\lim_{k\to\infty} u(x_k)H_0^{\frac{N-p}{p-1}}(x_k) = \CC \Hcap(\Omega)^{\frac{1}{p-1}} = \lim_{k\to\infty} u(y_k)H_0^{\frac{N-p}{p-1}}(y_k) ,
	$$
	i.e.
	$$
	\lim_{k\to\infty} r_k H_0^{\frac{N-p}{p-1}}(x_k) = \lim_{k\to\infty} r_k H_0^{\frac{N-p}{p-1}}(y_k) \,.
	$$
	By recalling the definition of $x_k$ and $y_k$ and \emph{Step 1}, we have
	$$
	\lim_{k\to\infty} H_0^{\frac{N-p}{p-1}}(x) =\lim_{k\to\infty} H_0^{\frac{N-p}{p-1}}(y),
	$$
	which implies that 
	$$
	H_0(x)=H_0(y)
	$$
	for every $x,y\in\bd U(\ls)$, i.e. $U(\ls)$ is Wulff shape.
	
	\medskip
	
	\emph{Conclusion:}	from \emph{Step 2.} we obtain that $U(r_k)$ is Wulff shape for any $k\ge 0$, which implies that the super level sets $U(s)$ are concentric Wulff shapes. In particular there exists $\beta >0$ such that $u=\beta H_0(x)^{1/\qq}$ for any $x \in \rn \setminus U(\ls)$. From Theorem \ref{thm_overdet} we conclude.

% % % % % % % % % % % % % % % % % % % % % % % % % % % %
%%%%%%%%%%%%%%%%%%%%%%%%%%%%%%%%%%%%%%%%%%%%%%%%%%%%%%%%%%%%

\section{Proof of Theorem \ref{main1}} \label{proof1}

	Let $\qq$ be given by \eqref{q}. Notice that for every $x_0\in\rn$ and every $R>0$, the concavity exponent of $\bhorx{x_0}{R}$ can be explicitly computed thanks to (\ref{potenziale-palla}) and it holds $\alpha(\bhorx{x_0}{R},p)= \qq$.
	
	We are going to prove that if the capacitary potential $u$ of the set $\Omega$ is $\qq$-concave then $\Omega$ is Wulff shape and this entails the desired result. Indeed if $u$ is $q$-concave, then $u$ is $s$-concave too, for every $s<q$.
	
	Assume that the function $u$ is $\qq$-concave. Since $\qq<0$, then $u^{\qq}$ is a convex function. 
	We denote by $V(t)$ the sublevel sets of the function $u^{\qq}$, i.e.  $V(t)=\{u^{\qq}\le t\}$; the superlevel sets of $u$ will be denoted by $U(t)$. Hence 
	$$
	U(t)=V(t^{\qq}) \,.
	$$ 
	Since $u^\qq$ is convex,  for every $t_0,t_1\in\R$ and every $\la\in[0,1]$ we have
	\begin{equation} \label{A}
	V((1-\la)t_0+\la t_1)\supseteq (1-\la)V(t_0)+\la V(t_1).
	\end{equation}
	Let $0<r<s<1$. By choosing $t_0=r^{\qq}, t_1=s^{\qq}$ and defining 
	\begin{equation} \label{t}
	t=((1-\la)r^{\qq}+\la s^{\qq})^{\frac 1{\qq}},
	\end{equation}
	\eqref{A} can be written as
	$$
	V(t^{\qq})\supseteq(1-\la)V(r^{\qq})+\la V(s^{\qq}),
	$$
	and hence
	$$
	U(t)\supseteq(1-\la)U(r)+\la U(s).
	$$
	From the monotonicity of the capacity and from Brunn-Minkowski inequality (\ref{BM}) 
	it follows
	\begin{eqnarray}\label{HcapOmegat}
	\Hcap(U(t)) &\ge& \Hcap((1-\la)U(r)+\la U(s))\nonumber\\
	&\ge& \Big((1-\la)\Hcap^{\frac 1{N-p}}(U(r)) + \la \Hcap^{\frac 1{N-p}}(U(s))\Big)^{N-p},
	\end{eqnarray}
	and, since  for every $r\in (0,1)$ 
	$$
	\Hcap(U(r))= {r^{1-p}}\Hcap(\Omega),
	$$
	inequality (\ref{HcapOmegat}) gives
	\begin{equation} \label{finale}
	\Hcap(\Omega) t^{1-p} \ge 
	\Hcap(\Omega) \Big( (1-\la)r^{\frac{1-p}{N-p}}+\la s^{\frac {1-p}{N-p}} \Big)^{N-p}.
	\end{equation}
	The definition of $t$ in \eqref{t} implies that the equality case holds in \eqref{finale} and this entails that the equality sign in the Brunn-Minkowski inequality (\ref{BM}) is attained. 
	Hence the superlevel set $U(r)$ is homothetic to $U(s)$ and Theorem \ref{thm-livelloomot} yields the conclusion.

%%%%%%%%%%%%%%%%%%%%%%%%%%%%%%%%%%%%%%%%%%%%%%%
%%%%%%%%%%%%%%%%%%%%%%%%%%%%%%%%%%%%%%%%%%%%%%%%%


\begin{thebibliography}{BSN2}
	
	%\bibitem{AF} E. Acerbi, N. Fusco, Regularity of minimizers of non-quadratic functionals: the case $1<p<2$, J. Math. Anal. Appl., {\bf 140} (1989), 115--135.
	%
	%\bibitem{AM} V. Agostiniani, L. Mazzieri, Riemannian aspects of potential theory, J. Math. Pures Appl., {\bf 104} (2015), 561--586.
	%
	\bibitem{AGHLV} M. Akman, J. Gong, J. Hineman, J. Lewis, A. Vogel, The Brunn-Minkowski inequality and a Minkowski problem for nonlinear capacity, preprint.
	%\bibitem{Al} A. D. Aleksandrov, Uniqueness theorems for surfaces
	%in the large V, Vestnik Leningrad Univ., {\bf 13}, no. 19 (1958), 5--8.
	%(English translation: Amer. Math. Soc. Translations, Ser. 2, {\bf 21} (1962),
	%412--415.)
	%
	%\bibitem{AKM} B. Avelin, T. Kuusi, G. Mingione, Nonlinear Calrder\'on-Zygmund theory in the limiting case. Preprint.
	%
	\bibitem{BL} H. J. Brascamp, E. H. Lieb, {\em On extensions
of the Brunn-Minkowski and Pr\'{e}kopa-Leindler theorems,
including inequalities for log-concave functions, and with an
application to the diffusion equation}, J. Funct. Anal. \textbf{22} (1976),
366-389.
	\bibitem{BLS} C. Bianchini, M. Longinetti, P. Salani, {\em Quasiconcave solutions to elliptic problems in convex rings}, Indiana Univ. Math. J. 58 no. 4 (2009), 1565Ð1589.
	%
	%\bibitem{BNP} G. Bellettini, M. Novaga, M. Paolini, On a crystalline variational problem, part I: first variation and global $L^\infty$ regularity, Arch. Ration. Mech. Anal., {\bf 157} (2001), 165--191.
	%
	%\bibitem{BP}  G. Bellettini, M. Paolini, Anisotropic motion by mean curvature in the context of Finsler geometry, Hokkaido Math. J., {\bf 25} (1996), 537--566.
	%
		
	\bibitem{BiCi} C. Bianchini, G. Ciraolo, Wulff shape characterization in overdetermined anisotropic problems, to appear on  Comm. Partial Differential Equations.
	
	%\bibitem{BiCiPali} C. Bianchini, G. Ciraolo, A note on an overdetermined problem for the capacitary potential. ``Geometric Properties for Parabolic and Elliptic PDE's'', Vol. 176 (2016) of the series Springer Proc. Math. Stat., 41--48. 
	%
	\bibitem{BCS} C. Bianchini, G. Ciraolo, P. Salani, An overdetermined problem for the anisotropic capacity, Calc. Var. Partial Differential Equations, 55:84 (2016).
	%
	%\bibitem{BNST} B. Brandolini, C. Nitsch, P. Salani, C. Trombetti, Serrin type overdetermined problems: an alternative proof, Arch. Rational Mech. Anal., {\bf 190} (2008), 267--280.
	%
	\bibitem{Bo} Ch. Borell, {\em Greenian potentials and concavity}, Math. Anal. 272 (1985), 155-160.
	\bibitem{CF} L.A. Caffarelli, A. Friedman, {\em Convexity of solutions of semilinear elliptic equations}, Duke Math. J. 52 (1985), 431Ð456. 
	%\bibitem{Ch} A. Chernov, Modern Crystallography III, Springer Ser. Solid-State Sci., vol. 36, Springer, Berlin, Heidelberg, 1984, softcover reprint of the original 1st edition.
	%
	%\bibitem{CS} A. Cianchi, P. Salani, Overdetermined anisotropic elliptic problems, Math. Ann., {\bf 345} (2009), 859--881.
	%
	%\bibitem{CMS} G. Ciraolo, R. Magnanini, S. Sakaguchi, Symmetry of solutions of elliptic and parabolic equations with a level surface parallel to the boundary, J. Eur. Math. Soc., {\bf 17} (2015), 2789--2804. 
	%
	%\bibitem{CV} G. Ciraolo, L. Vezzoni, A rigidity problem on the round sphere. To appear in Comm. Contemp. Math. 
	%
	\bibitem{CS} A. Colesanti and P. Salani, {\em The
Brunn-Minkowski inequality for $p$-capacity of convex bodies}, Math. Ann. \textbf{327} (2003), 459-479.

	%\bibitem{CFV}  M. Cozzi, A. Farina, E. Valdinoci, Gradient bounds and rigidity results for singular, degenerate, anisotropic partial differential equations, Comm. Math. Phys., {\bf 331} (2014), 189--214.
	%
	%\bibitem{CFV2} M. Cozzi, A. Farina, E. Valdinoci, Monotonicity formulae and classification results for singular, degenerate, anisotropic PDEs, Adv. Math., {\bf 293} (2016), 343--381.
	%
	
	%\bibitem{CIL} M.G. Crandall, H. Ishii, P.L. Lions, User's guide to viscosity solutions of second order partial differential equations, 1992 **************
	%%\bibitem{GdP1}  F. Della Pietra, N. Gavitone, Anisotropic elliptic problems involving Hardy-type potentials. J. Math. Anal. Appl. 397 (2013), no. 2, 800-813.
	%%
	%%\bibitem{GdP2}  F. Della Pietra, N. Gavitone, Anisotropic elliptic equations with general growth in the gradient and Hardy-type potentials. J. Differential Equations 255 (2013), no. 11, 3788-3810.
	%%
	%% \bibitem{GdP3}  F. Della Pietra, N. Gavitone, Sharp bounds for the first eigenvalue and the torsional rigidity related to some anisotropic operators. Math. Nachr. 287 (2014), no. 2-3, 194-209.
	
	
	%\bibitem{DiB} E. Dibenedetto, $C^{1,\alpha}$-local regularity of weak solutions of degenerate elliptic equations, Nonlinear Anal., {\bf 7} (1983), 827--859.
	%
	%\bibitem{EM} L. Esposito, G. Mingione, Some remarks on the regularity of weak solutions of degenerate elliptic systems, Rev. Mat. Complutense, {\bf 11} (1998), 203--219.
	%
	%\bibitem{FK} A. Farina, B. Kawohl, Remarks on an overdetermined boundary value problem, Calculus of Variations and Partial Differential Equations, 31:351-357 (2008).
	%%\bibitem{FV} A. Farina, E. Valdinoci, Gradient bounds for anisotropic partial differential equations, Calc. Var. Partial Differential Equations, 49 (2014), no. 3-4, 923--936.
	
	\bibitem{FK} V. Ferone, B. Kawohl, Remarks on a Finsler-Laplacian, Proc. Am. Math. Soc. 137 (2009), 247-253.
	
	%
	%\bibitem{FC} A. Figalli, M. Colombo, An excess-decay result for a class of degenerate elliptic equations, Discrete Contin. Dyn. Syst. Ser. S, {\bf 7} (2014), 631--652.
	%
	%\bibitem{FGK} I. Fragal\`a, F. Gazzola, B. Kawohl, Overdetermined problems with possibly degenerate ellipticity, a geometric approach, 
	%Math. Z., {\bf 254} (2006), 117-.132.
	%
	%\bibitem{GL} N. Garofalo, N. Lewis, A symmetry result related to some overdetermined boundary value problems, Am. J. Math. 111 (1989), 9-33.
	\bibitem{G}M. Gabriel, {\em A result concerning convex level--surfaces 
of three--dimensional harmonic functions}, London Math. Soc. J. \textbf{32} (1957),
286--294.
	
	\bibitem{GL1} N. Garofalo and F-.H. Lin, \emph{Monotonicity properties of variational integrals, $A_p$ weights and unique continuation,} Indiana Univ. Math. J., 35 (1986), 245-268.
	
	\bibitem{GL2} N. Garofalo, N. and F.-H. Lin, \emph{Unique continuation for elliptic operators: a geometric-variational approach}, Comm. Pure Appl. Math., 40 (1987), 347-366.
	
	\bibitem{HNST} A. Henrot, C. Nitsch, C. Trombetti, P. Salani, {Optimal concavity of the torsion function}, to appear in J. Optim. Theory Appl.
	
	
\bibitem{Ke} A. U. Kennington, \emph{Power concavity and boundary value problems}, Indiana Univ. Math. J.
34 (1985), 687-704.

\bibitem{L} J. Lewis,
\emph{Capacitary functions in convex rings},
Arch. Rational Mech. Anal. {\bf 66} (1977), 201--224.
		
	%
	%
	%
	%\bibitem{GT} D. Gilbarg, N.S. Trudinger, Elliptic Partial Differential Equations of Second Order, Second Edition, Springer 1997.
	%
	%\bibitem{Gi} E. Giusti, Metodi diretti nel calcolo delle variazioni. (Italian) [Direct methods in the calculus of variations], Unione Matematica Italiana, Bologna, 1994.
	%
	%\bibitem{HS} A. Henrot, H. Shahgholian, Existence of classical solutions to a free boundary problem for the $p$-Laplace operator: (II) the interior case, Indiana Univ. Nath. J., 49 (2000), 311-323.
	%
	%%\bibitem{HeLi} Y. J. He, H. Z. Li, Integral formula of Minkowski type and new characterization of the Wulff shape, Acta Math. Sin. (Engl. Ser.) 24 (2008), n. 4, 697-704.
	%
	%\bibitem{HLMG} Y. He, H. Li, H. Ma, J. Ge, Compact embedded hypersurfaces with constant higher order anisotropic mean curvatures, Indiana Univ. Math. J. 58 No. 2 (2009), 853-868.
	%
	%\bibitem{Ke} S. Kesavan, Symmetrization \& applications, Series in Analysis, 3, World Scientific
	%Publishing Co. Pte. Ltd., Hackensack, NJ, 2006.
	%
	%\bibitem{KP} S. Kumaresan, J. Prajapat, Serrin's result for hyperbolic space and sphere, Duke Math. J., 91 (1998), 17-28.
	%
	%\bibitem{LU} O. A. Ladyzhenskaya, N. N. Uraltseva, Linear and Quasilinear Elliptic Equations, Accademic Press, New York, 1968.
	%
	%\bibitem{Le} J. L. Lewis, Capacitary functions in convex rings, Arch. Rational Mech. Anal., 66 (1977), no. 3, 201-224.
	%
	%\bibitem{LL} E.H. Lieb, M. Loss, Analysis, Graduate Studies in Mathematics, 14. American Mathematical Society, Providence, RI, (1997).
	%
	%\bibitem{Lieberman} G. M. Lieberman, Boundary regularity for solutions of degenerate elliptic equations, Nonlinear Anal., {\bf 12} (1988), 1203--1219.
	%
	%\bibitem{LLZ} B. Lv, F. Li, W. Zou,  Overdetermined boundary value problems with strongly nonlinear elliptic PDE, Electron. J. Qual. Theory Differ. Equ. 2012, No. 10, 17 pp.
	\bibitem{ML} L.G. Makar-Limanov, {\em The solution of the Dirichlet problem for the equation $\Delta u=-1$ in a convex region}, Mat. Zametki 9 (1971) 89-92 (Russian).
English translation in Math. Notes 9 (1971), 52-53.
	%
	%\bibitem{Molz} R. Molzon, Symmetry and overdetermined boundary value problems, Forum Math., 3 (1991), 143-156.
	%
	%\bibitem{Mo} J. Mossino, Inegalit\'es Isop\'erim\'etriques et Applications en Physique, Hermann (1984).
	%
	%\bibitem{NP} M. Novaga, E. Paolini, A computational approach to fractures in crystal growth, Atti Accad.
	%Naz. Lincei Cl. Sci. Fis. Mat. Natur., 10 (1999), 47-56.
	%
	%\bibitem{OBGXY} S. Osher, M. Burger, D. Goldfarb, J. Xu, W. Yin, An iterative regularization method for
	%total variation-based image restoration, Multiscale Model. Simul., 4 (2005), 460-489.
	%
	%\bibitem{PS} G. P\'olya,  G. Szeg\"o, Isoperimetric inequalities in Mathematical Physics, Ann. of Math. Studies, 27,
	%Princeton University Press, Princeton, 1951.
	%
	%\bibitem{Re} W. Reichel, Radial Symmetry for Elliptic Boundary-Value Problems on Exterior Domains, Arch. Rational Mech. Anal., 137 (1997), no. 4, 381-394.
	%
	%
	\bibitem{Sa} P. Salani, {\em A characterization of balls through optimal concavity for potential functions}, Proc. AMS 143 (1) (2015), 173--183.	
	%\bibitem{Sc} R. Schneider, Convex Bodies: The Brunn-Minkowski Theory, Cambridge University Press, Cambridge (1993).
	%
	%\bibitem{Se_acta}  J. Serrin, Local behavior of solutions of quasi-linear equations, Acta Math., 111 (1964), 247-302.
	%
	%\bibitem{Se} J. Serrin, A symmetry problem in potential theory, Arch. Rational Mech. Anal., 43 (1971), 304-318.
	%
	%\bibitem{Sz} G. Szeg\'o, \"Uber einige Extremalaufgaben der Potentialtheorie, Math. Z., 31 (1930), 583-593.
	%
	%\bibitem{Tay} J. Taylor, Crystalline variational problems, Bull. Amer. Math. Soc., 84 (1978) 568-588.
	%
	%\bibitem{TCH} J.E. Taylor, J.W. Cahn, C.A Handwerker, Geometric models of crystal growth. Acta Metall., 40 (1992), 1443-1474.
	%
	%%
	%\bibitem{To} P. Tolksdorf, Regularity for a more general class of quasilinear elliptic equations, J. Diff. Eq., 51 (1984), 126-150.
	%%
	%%
	%
	%\bibitem{WX} G. Wang, C. Xia, A characterization of the Wulff shape by an overdetermined anisotropic PDE, Arch. Ration. Mech. Anal. 199 (2011), no. 1, 99-115.
	%%
	%%\bibitem{WX1} G. Wang, C. Xia, A characterization of the Wulff shape by an overdetermined anisotropic PDE. Arch. Ration. Mech. Anal. 199 (2011), no. 1, 99-115.
	%%
	%%\bibitem{WX2}  G. Wang, C. Xia, An optimal anisotropic Poincar\'e inequality for convex domains. Pacific J. Math. 258 (2012), no. 2, 305-325.
	%%
	%
	%\bibitem{We} H. Weinberger, Remark on the preceding paper of Serrin, Arch. Ration. Mech. Anal. 43 (1971), 319-320.
	%
	%\bibitem{Wu} G. Wulff, Zur Frage der Geschwindigkeit
	%des Wachstums und der Aufl\"osung der Kristallfl\"aschen, Z.
	%Krist., 34 (1901), 449--530.
	%%
	%%%\bibitem{XiaThesis} C. Xia, On a class of anisotropic problems, Dissertation zur Erlangung des Doktorgrades der Fakult\"at Mathematik und Physik der Albert-Ludwigs-Universit\"at Freiburg im Breisgau, April 2012.
	
	
\end{thebibliography}
\end{document}